\documentclass[twosided,reqno]{amsart}

\usepackage{amsmath}
\usepackage{amssymb}
\usepackage{amsthm}

\newcommand{\e}{e_{\lambda}}

\theoremstyle{theorem}
\newtheorem{theorem}{\scshape Theorem }[section]

\theoremstyle{definition}

\numberwithin{equation}{section}

\begin{document}

\title[]{A note on multi-Euler-Genocchi and degenerate multi-Euler-Genocchi polynomials}

\author{Taekyun  Kim$^{1}$}
\address{$^{1}$ Department of Mathematics, Kwangwoon University, Seoul 139-701, Republic of Korea}
\email{tkkim@kw.ac.kr}

\author{Dae San Kim$^{2}$}
\address{$^{2}$ Department of Mathematics, Sogang University, Seoul 121-742, Republic of Korea}
\email{dskim@sogang.ac.kr}

\author{Jin-Woo Park$^{3,*}$}
\address{Department of Mathematics Education, Daegu University,
38453, Republic of Korea.}
\email{a0417001@daegu.ac.kr}
\thanks{$^*$ Corresponding author}

\author{Jongkyum Kwon$^{4}$}
\address{$^{4}$ Department of Mathematics Education and ERI, Gyeongsang National University, Jinju, Gyeongsangnamdo, 52828, Republic of Korea}
\email{mathkjk26@gnu.ac.kr}



\begin{abstract}
Recently, introduced are the generalized Euler-Genocchi and generalized degenerate Euler-Genocchi polynomials. The aim of this note is to study the multi-Euler-Genocchi and degenerate multi-Euler-Genocchi polynomials which are defined by means of the multiple logarithm and generalize respectively the generalized Euler-Genocchi and generalized degenerate Euler-Genocchi polynomials. Especially, we express the former by the generalized Euler-Genocchi polynomials, the multi-Stirling numbers of the first kind and Stirling numbers of the second kind, and the latter by the generalized degenerate Euler-Genocchi polynomials, the multi-Stirling numbers of the first kind and Stirling numbers of the second kind.
\end{abstract}

\maketitle

\markboth{\centerline{\scriptsize multi-Euler-Genocchi and degenerate multi-Euler-Genocchi polynomials}}
{\centerline{\scriptsize T. Kim, D. S. Kim, J.-W. Park, and J. Kwon}}

\section{Introduction}
Carlitz began a study of degenerate versions of Bernoulli and Euler polynomials, namely the degenerate Bernoulli and degenerate Euler polynomials (see \cite{02,03}). Recently, explorations for various degenerate versions of some special numbers and polynomials have regained interests of many mathematicians and yielded lots of fascinating and fruitful results. Indeed, this quest for degenerate versions even led to the development of degenerate umbral calculus \cite{07}, degenerate $q$-umbral calculus \cite{11} and degenerate gamma function \cite{09}. \par
Goubi introduced the generalized Euler-Genocchi polynomials (of order $\alpha$) in \cite{05}. A degenerate version of those polynomials, namely the generalized degenerate Euler-Genocchi polynomials, is investigated in \cite{10}. Here we study a`multi-version' of them, namely the multi-Euler-Genocchi and degenerate multi-Euler-Genocchi polynomials. The generating function of the Stirling numbers of the first kind is given by the usual logarithm (see \eqref{9}). Naturally, the multi-Stirling numbers of the first kind (see \eqref{13}) are defined by means of the multiple logarithm (see \eqref{7}). In the same way, both the multi-Euler-Genocchi and the degenerate multi-Euler-Genocchi polynomials are defined by using the multiple logarithm, and the latter is a degenerate version of the former  (see \eqref{14}, \eqref{21}). \par
 
The aim of this note is to study the multi-Euler-Genocchi and degenerate multi-Euler-Genocchi polynomials which generalize respectively the generalized Euler-Genocchi (see \eqref{3}) and generalized degenerate Euler-Genocchi polynomials (see \eqref{6}). \par
The outline of this note is as follows. In Theorem 2.1, the value of the multi-Euler-Genocchi polynomial at 1 is expressed in terms of the multi-Stirling numbers of the first kind and the Stirling numbers of the second kind. The multi-Euler-Genocchi polynomials are represented in terms of the generalized Euler-Genocchi polynomials, the multi-Stirling numbers of the first kind and the Stirling numbers of the second kind in Theorem 2.2. Likewise, the degenerate multi-Euler-Genocchi polynomials are expressed in terms of the generalized degenerate Euler-Genocchi polynomials, the multi-Stirling numbers of the first kind and the Stirling numbers of the second kind in Theorem 2.4. A distribution type formula is derived for the degenerate multi-Euler-Genocchi polynomials in Theorem 2.6. In the rest of this section, we recall the facts that are needed throughout this paper. \par

The {\it{Euler polynomials}} are defined by 
\begin{equation}\label{1}
\frac{2}{e^t+1}e^{xt}=\sum_{n=0} ^{\infty} E_n(x)\frac{t^n}{n!},{\text{ (see \cite{01, 04, 14})}}.
\end{equation}
When $x=0$, $E_n=E_n(0)$, $(n\geq0)$ are called the {\it{Euler numbers}}. The {\it{Genocchi polynomials}} are given by 
\begin{equation}\label{2}
\frac{2t}{e^t+1}e^{xt}=\sum_{n=0} ^{\infty} G_n(x)\frac{t^n}{n!},{\text{ (see \cite{01, 04, 14})}}.
\end{equation}
When $x=0$, $G_n=G_n(0)$ are called the Genocchi numbers. We note that $G_n \in {\mathbb{Z}}$, $(n\geq0)$. \\
Recently, the {\it{generalized Euler-Genocchi polynomials}} are introduced as
\begin{equation}\label{3}
\frac{2t^r}{e^t+1}e^{xt}=\sum_{n=0} ^{\infty} A_n ^{(r)}(x)\frac{t^n}{n!},{\text{ (see \cite{05})}},
\end{equation}
where $r$ is a nonnegative integer. Note that $A_n ^{(0)}(x)=E_n(x)$, $A_n ^{(1)}(x)=G_n(x)$, $(n\geq0)$. \par
For any nonzero $\lambda\in{\mathbb{R}}$, the {\it{degenerate exponentials}} are defined by
\begin{equation}\label{4}
\e ^x(t)=\left(1+\lambda t\right)^{\frac{x}{\lambda}}=\sum_{n=0} ^{\infty} (x)_{n,\lambda}\frac{t^n}{n!},{\text{ (see \cite{08,12})}},
\end{equation}
where $(x)_{0,\lambda}=1$, $(x)_{n,\lambda}=x(x-\lambda)\cdots(x-(n-1)\lambda)$, $(n\geq1)$. In particular, for $x=1$, $\e(t)=\e ^1(t)= \sum_{n=0} ^{\infty} \frac{(1)_{n,\lambda}}{n!}t^n$. \par
In \cite{02,03}, Carlitz introduced the {\it{degenerate Euler polynomials}} given by 
\begin{equation}\label{5}
\frac{2}{\e(t)+1}\e ^x(t)=\sum_{n=0} ^{\infty} {\mathcal{E}}_{n,\lambda}(x)\frac{t^n}{n!}.
\end{equation}
Note that $\lim_{n\rightarrow 0} {\mathcal{E}}_{n,\lambda}(x)=E_n(x)$, $(n\geq0)$.
In view of \eqref{2}, the degenerate Genocchi polynomials are defined by
\begin{equation*}
\frac{2t}{\e(t)+1}\e ^x(t)=\sum_{n=0} ^{\infty} G_{n,\lambda}(x) \frac{t^n}{n!},{\text{ (see \cite{10})}}.
\end{equation*}
Recently, the {\it{generalized degenerate Euler-Genocchi polynomials}} are defined by
\begin{equation}\label{6}
\frac{2t^r}{\e(t)+1}\e ^x (t)=\sum_{n=0} ^{\infty} A_{n,\lambda} ^{(r)}(x)\frac{t^n}{n!},{\text{ (see \cite{10})}}.
\end{equation}
When $x=0$, $A_{n,\lambda} ^{(r)}=A_{n,\lambda} ^{(r)}(0)$ are called the {\it{generalized degenerate Euler-Genocchi numbers}}. Note that $\lim_{\lambda\rightarrow 0} A_{n,\lambda} ^{(r)}(x)=A_n ^{(r)}(x)$, $(n\geq0)$. \par
For $k_1,k_2,\ldots,k_r\in{\mathbb{Z}}$, the {\it{multiple logarithm}} is defined by
\begin{equation}\label{7}
Li_{k_1,k_2,\ldots,k_r}(t)=\sum_{0<n_1<n_2<\cdots<n_r} \frac{t^{n_r}}{n_1 ^{k_1} n_2 ^{k_2} \cdots n_r ^{k_r}},~(|t|<1),{\text{ (see \cite{06})}}.
\end{equation}
The {\it{multi-Bernoulli numbers}} are defined by
\begin{equation}\label{8}
\frac{r!Li_{k_1,k_2,\ldots,k_r}\left(1-e^{-t}\right)}{(\e(t)-1)^r}=\sum_{n=0} ^{\infty} B_n ^{(k_1,k_2,\ldots,k_r)}\frac{t^n}{n!},{\text{ (see \cite{06,13})}}.
\end{equation}
It is well known that the {\it{Stirling numbers of the first kind}} are given by 
\begin{equation}\label{9}
\frac{1}{k!}\left(\log(1+t)\right)^k=\sum_{n=k} ^{\infty} S_1(n,k) \frac{t^n}{n!},~(k\geq0),{\text{ (see \cite{04, 14})}}.
\end{equation}
As the inversion formula of \eqref{9}, the {\it{Stirling numbers of the second kind}} are defined by 
\begin{equation}\label{10}
\frac{1}{k!}\left(e^t-1\right)^k =\sum_{n=k} ^{\infty} S_2(n,k)\frac{t^n}{n!},~(k\geq0),{\text{ (see \cite{14})}}.
\end{equation} \par
From \eqref{7}, we note that
\begin{equation}\label{11}
\begin{split}
\frac{d}{dt}Li_{k_1,k_2,\ldots,k_r}(t)=&\frac{d}{dt}\sum_{0<n_1<n_2<\cdots<n_r} \frac{t^n_r}{n_1 ^{k_1} n_2 ^{k_2}\cdots n_r ^{k_r}}\\
=&\frac{1}{t}Li_{k_1,\ldots,k_{r-1},k_r-1}(t),{\text{ (see \cite{08})}}.
\end{split}
\end{equation}
By \eqref{11}, we get
\begin{equation}\label{12}
Li_{\underbrace{1,1,\ldots,1}_{r-times}}(t)=\frac{1}{r!}\left(-\log(1-t)\right)^r=\sum_{n=r} ^{\infty} (-1)^{n-r}S_1(n,r)\frac{t^n}{n!}.
\end{equation} \par
In light of \eqref{9} and \eqref{12}, the {\it{multi-Stirling numbers of the first kind}} are defined by
\begin{equation}\label{13}
Li_{k_1,k_2,\ldots,k_r}(t)=\sum_{n=r} ^{\infty} S_1 ^{(k_1,k_2,\ldots,k_r)}(n,r)\frac{t^n}{n!},{\text{ (see \cite{13})}}.
\end{equation}
Note from \eqref{12} and \eqref{13} that 
\begin{equation*}
S_1 ^{\overbrace{(1,1,\ldots,1)}^{r-times}}(n,r)=(-1)^{n-r}S_1(n,r), \,\,(n,r\geq0).
\end{equation*}

\section{Multi-Euler-Genocchi and degenerate multi-Euler-Genocchi polynomials}

We consider the {\it{multi-Euler-Genocchi polynomials}} given by
\begin{equation}\label{14}
\frac{2r!}{e^t+1}Li_{k_1,k_2,\ldots,k_r}\left(1-e^{-t}\right)e^{xt}=\sum_{n=0} ^{\infty} A_n ^{(k_1,k_2,\ldots,k_r)}(x)\frac{t^n}{n!},
\end{equation}
where $r$ is a nonnegative integer. When $x=0$, $A_n ^{(k_1,k_2,\ldots,k_r)}=A_n ^{(k_1,k_2,\ldots,k_r)}(0)$ are called the {\it{multi-Euler-Genocchi numbers}}.

From \eqref{12} and \eqref{14}, we note that
\begin{equation}\label{15}
\begin{split}
\sum_{n=0} ^{\infty} A_n ^{\overbrace{(1,1,\ldots,1)}^{r-times}}(x)\frac{t^n}{n!}=&\frac{2r!}{e^t+1}Li_{\underbrace{1,1,\ldots,1}_{r-times}}\left(1-e^{-t}\right)e^{xt}\\
=&\frac{2t^r}{e^t+1}e^{xt}=\sum_{n=0} ^{\infty} A_n ^{(r)}(x)\frac{t^n}{n!}.
\end{split}
\end{equation}
Thus, by \eqref{15}, we get
\begin{equation}\label{16}
A_n ^{(\overbrace{1,1,\ldots,1}^{r-times})}(x) =A_n ^{(r)}(x),~(n\geq0).
\end{equation}

From \eqref{13} and \eqref{15}, we have
\begin{equation}\label{17}
\begin{split}
&\sum_{n=0} ^{\infty} A_n ^{(k_1,k_2,\ldots,k_r)}(1)\frac{t^n}{n!}=\frac{r!}{1+\frac{1}{2}\left(e^{-t}-1\right)}Li_{k_1,k_2,\ldots,k_r}\left(1-e^{-t}\right)\\
=&r!\sum_{l=0} ^{\infty}\left(\frac{1}{2}\right)^l(-1)^l \left(e^{-t}-1\right)^l\sum_{m=r} ^{\infty} S_1 ^{(k_1,k_2,\ldots,k_r)}(m,r)(-1)^m\frac{1}{m!}\left(e^{-t}-1\right)^m\\
=&r!\sum_{l=0}^{\infty}\sum_{m=r}^{\infty}\left(\frac{1}{2}\right)^{l}S_{1}^{(k_1,k_2,\dots,k_r)}(m,r)(-1)^{l+m}\frac{(l+m)!}{m!}\frac{1}{(l+m)!}(e^{-t}-1)^{l+m}\\
=&r!\sum_{l=0}^{\infty}\sum_{m=r}^{\infty}\left(\frac{1}{2}\right)^{l}S_{1}^{(k_1,k_2,\dots,k_r)}(m,r)(-1)^{l+m}\frac{(l+m)!}{m!}\sum_{n=l+m}^{\infty}S_{2}(n,l+m)(-1)^{n}\frac{t^n}{n!}\\
=&\sum_{n=r} ^{\infty} \left(r!\sum_{k=r} ^n \sum_{m=r} ^k \left(\frac{1}{2}\right)^{k-m}(-1)^{n-k}S_1 ^{\left(k_1,k_2,\ldots,k_r\right)}(m,r)\frac{k!}{m!}S_2(n,k)\right)\frac{t^n}{n!}.
\end{split}
\end{equation}
Therefore, by comparing the coefficient on both sides of \eqref{17}, we obtain the following theorem.
\begin{theorem}
For $n,r\in{\mathbb{Z}}$ with $n \geq r \geq 0$, we have
\begin{equation*}
A_n ^{(k_1,k_2,\ldots,k_r)}(1)=r!\sum_{k=r} ^n \sum_{m=r} ^k \left(\frac{1}{2}\right)^{k-m}(-1)^{n-k}\frac{k!}{m!}S_1 ^{(k_1,k_2,\ldots,k_r)}(m,r)S_2(n,k),
\end{equation*}
and, for $0 \le n<r$, we have
\begin{equation*}
A_n ^{(k_1,k_2,\ldots,k_r)}(1)=0.
\end{equation*}
\end{theorem}

By \eqref{14}, we get
\begin{equation}\label{18}
\begin{split}
&\sum_{n=0} ^{\infty} A_n ^{(k_1,k_2,\ldots,k_r)}(x)\frac{t^n}{n!}=\frac{2r!}{e^t+1}Li_{k_1,k_2,\ldots,k_r}\left(1-e^{-t}\right)e^{xt}\\
=&r!\frac{2}{e^t+1}e^{xt}\sum_{l=r} ^{\infty} S_1 ^{(k_1,k_2,\ldots,k_r)}(l,r)(-1)^l\frac{1}{l!}\left(e^{-t}-1\right)^l\\
=&\frac{r!2e^{xt}}{e^t+1}\sum_{l=r} ^{\infty} S_1 ^{(k_1,k_2,\ldots,k_r)}(l,r)(-1)^l\sum_{m=l} ^{\infty} S_2(m,l)(-1)^m\frac{t^m}{m!}\\
=&\frac{r!2e^{xt}}{e^t+1}\sum_{m=r}^{\infty}\left(\sum_{l=r} ^m S_1^{(k_1,k_2,\ldots,k_r)}(l,r)S_2(m,l)(-1)^{m-l}\right)\frac{t^m}{m!}\\
=&\frac{r!2t^re^{xt}}{e^t+1}\sum_{m=0} ^{\infty}\left(\sum_{l=r} ^{m+r}S_1 ^{(k_1,k_2,\ldots,k_r)}(l,r)\frac{S_2(m+r,l)m!}{(m+r)!}(-1)^{m+r-l}\right)\frac{t^m}{m!}\\
=&\sum_{j=0} ^{\infty} A_j ^{(r)}(x)\frac{t^j}{j!}\sum_{m=0} ^{\infty}\left(\sum_{l=r} ^{m+r}\frac{S_1 ^{(k_1,k_2,\ldots,k_r)}(l,r)}{\binom{m+r}{r}}S_2(m+r,l)(-1)^{m+r-l}\right)\frac{t^m}{m!}\\
=&\sum_{n=0} ^{\infty}\left(\sum_{m=0} ^n \binom{n}{m}A_{n-m} ^{(r)}(x)\sum_{l=r} ^{m+r} \frac{S_1 ^{(k_1,k_2,\ldots,k_r)}(l,r)}{\binom{m+r}{r}}S_2(m+r,l)(-1)^{m+r-l}\right)\frac{t^n}{n!}.
\end{split}
\end{equation}
Therefore, by comparing the coefficients on both sides of \eqref{18}, we obtain the following theorem.
\begin{theorem}
For $n \geq0$, we have
\begin{equation*}
A_n ^{(k_1,k_2,\ldots,k_r)}(x)=\sum_{m=0} ^n \binom{n}{m}A_{n-m} ^{(r)}(x)\sum_{l=r} ^{m+r} \frac{S_1 ^{(k_1,k_2,\ldots,k_r)}(l,r)}{\binom{m+r}{r}}S_2(m+r,l)(-1)^{m+r-l}.
\end{equation*}
\end{theorem}

From \eqref{13}, we note that
\begin{equation}\label{19}
\begin{split}
&Li_{k_1,k_2,\ldots,k_r}\left(1-e^{-t}\right)=\sum_{k=r} ^{\infty} S_1 ^{(k_1,k_2,\ldots,k_r)}(k,r)\frac{1}{k!}\left(1-e^{-t}\right)^k\\
=&\sum_{k=r} ^{\infty} S_1 ^{(k_1,k_2,\ldots,k_r)}(k,r)\sum_{n=k} ^{\infty} (-1)^{n-k}S_2(n,k)\frac{t^n}{n!}\\
=&\sum_{n=r} ^{\infty}\left(\sum_{k=r} ^n S_1 ^{(k_1,k_2,\ldots,k_r)}(k,r)(-1)^{n-k}S_2(n,k)\right)\frac{t^n}{n!}.
\end{split}
\end{equation}
On the other hand, by \eqref{14}, we get
\begin{equation}\label{20}
\begin{split}
Li_{k_1,k_2,\ldots,k_r}\left(1-e^{-t}\right)=&\frac{1}{2r!}\frac{2r!}{e^t+1}Li_{k_1,k_2,\ldots,k_r}\left(1-e^{-t}\right)\left(e^t+1\right)\\
=&\frac{1}{2r!}\sum_{n=0} ^{\infty}\left(A_n ^{(k_1,k_2,\ldots,k_r)}(1)+A_n ^{(k_1,k_2,\ldots,k_r)}\right)\frac{t^n}{n!}.
\end{split}
\end{equation}
Therefore, by \eqref{19} and \eqref{20}, we obtain the following theorem.
\begin{theorem}
For $n,r\geq0$, we have
\begin{equation*}
A_n ^{(k_1,k_2,\ldots,k_r)}(1)+A_n ^{(k_1,k_2,\ldots,k_r)}=
\begin{cases}
2r!\sum_{k=r} ^n S_1 ^{(k_1,k_2,\ldots,k_r)}(k,r)(-1)^{n-k}S_2(n,k), & {\text{if }}n\geq r,\\
0,& {\text{if }}0 \le n<r.
\end{cases}
\end{equation*}
\end{theorem}

Now, we consider the {\it{degenerate multi-Euler-Genocchi polynomials}} given by
\begin{equation}\label{21}
\frac{2r!}{\e(t)+1}Li_{k_1,k_2,\ldots,k_r}\left(1-e^{-t}\right)\e^x(t)=\sum_{n=0} ^{\infty} A_{n,\lambda} ^{(k_1,k_2,\ldots,k_r)}(x)\frac{t^n}{n!}.
\end{equation}
When $x=0$, $A_{n,\lambda} ^{k_1,k_2,\ldots,k_r)}=A_{n,\lambda} ^{(k_1,k_2,\ldots,k_r)}(0)$ are called the {\it{degenerate multi-Euler-Genocchi numbers}}.
Thus, by \eqref{21}, we get
\begin{equation}\label{22}
\begin{split}
\sum_{n=0} ^{\infty} A_{n,\lambda} ^{(\overbrace{1,1,\ldots,1}^{r-times})}(x)\frac{t^n}{n!}=&\frac{2r!}{\e(t)+1}Li_{\underbrace{1,1,\ldots,1}_{r-times}}\left(1-e^{-t}\right)\e^x(t)\\
=&\frac{2t^r}{\e(t)+1}\e^x(t)=\sum_{n=0} ^{\infty} A_{n,\lambda} ^{(r)}(x)\frac{t^n}{n!}.
\end{split}
\end{equation}
From \eqref{22}, we have
\begin{equation}\label{23}
A_{n,\lambda} ^{(\overbrace{1,1,\ldots,1}^{r-times})}(x)=A_{n,\lambda} ^{(r)}(x),~(n\geq0).
\end{equation}

By \eqref{21}, we get
\begin{equation}\label{24}
\begin{split}
&\sum_{n=0} ^{\infty} A_{n,\lambda} ^{(k_1,k_2,\ldots,k_r)}(x)\frac{t^n}{n!}=\frac{2r!}{\e(t)+1}Li_{k_1,k_2,\ldots,k_r}\left(1-e^{-t}\right)\e^x(t)\\
=&\frac{2r!}{\e(t)+1}\e^x(t)\sum_{m=r} ^{\infty} S_1 ^{(k_1,k_2,\ldots,k_r)}(m,r)\frac{1}{m!}\left(1-e^{-t}\right)^m\\
=&\frac{2r!}{\e(t)+1}\e^x(t)\sum_{m=r} ^{\infty} S_1 ^{(k_1,k_2,\ldots,k_r)}(m,r)\sum_{l=m} ^{\infty} S_2(l,m)(-1)^{l-m}\frac{t^l}{l!}\\
=&\frac{2r!}{\e(t)+1}\e^x(t)\sum_{l=r} ^{\infty}\left(\sum_{m=r} ^l (-1)^{l-m}S_2(l,m)S_1^{(k_1,k_2,\ldots,k_r)}(m,r)\right)\frac{t^l}{l!}\\
=&\sum_{j=0} ^{\infty} A_{j,\lambda} ^{(r)}(x)\frac{t^j}{j!}\sum_{l=0} ^{\infty}\left(\sum_{m=r} ^{l+r}\frac{S_2(l+r,m)}{\binom{l+r}{r}}(-1)^{l-m-r}S_1 ^{(k_1,k_2,\ldots,k_r)}(m,r)\right)\frac{t^l}{l!}\\
=&\sum_{n=0} ^{\infty}\left(\sum_{l=0} ^n \binom{n}{l}A_{n-l,\lambda} ^{(r)}(x)\sum_{m=r} ^{l+r}\frac{S_2(l+r,m)}{\binom{l+r}{l}}(-1)^{l-m-r}S_1 ^{(k_1,k_2,\ldots,k_r)}(m,r)\right)\frac{t^n}{n!}.
\end{split}
\end{equation}
Therefore, by comparing the coefficients on both sides of \eqref{24}, we obtain the following theorem.
\begin{theorem}
for $n \geq 0$, we have
\begin{equation*}
A_{n,\lambda} ^{(k_1,k_2,\ldots,k_r)}(x)=\sum_{l=0} ^n \binom{n}{l}A_{n-l,\lambda} ^{(r)}(x)\sum_{m=r} ^{l+r}\frac{S_2(l+r,m)}{\binom{l+r}{l}}(-1)^{l-m-r}S_1 ^{(k_1,k_2,\ldots,k_r)}(m,r).
\end{equation*}
\end{theorem}

From \eqref{21}, we note that
\begin{equation}\label{25}
\begin{split}
&\sum_{n=0} ^{\infty} A_{n,\lambda} ^{(k_1,k_2,\ldots,k_r)}(x)\frac{t^n}{n!}=\frac{2r!Li_{k_1,k_2,\ldots,k_r}\left(1-e^{-t}\right)}{\e(t)+1}\e^x(t)\\
=&\sum_{k=0} ^{\infty} A_{k,\lambda} ^{(k_1,k_2,\ldots,k_r)}\frac{t^k}{k!}\sum_{m=0} ^{\infty}(x)_{m,\lambda}\frac{t^m}{m!}\\
=&\sum_{n=0} ^{\infty} \left(\sum_{k=0} ^n \binom{n}{k} A_{k,\lambda} ^{(k_1,k_2,\ldots,k_r)}(x)_{n-k,\lambda}\right)\frac{t^n}{n!}.
\end{split}
\end{equation}
Thus, by \eqref{25}, we obtain the following theorem.
\begin{theorem}
For $n\geq0$, we have
\begin{equation*}
A_{n,\lambda} ^{(k_1,k_2,\ldots,k_r)}(x)=\sum_{k=0} ^n \binom{n}{k} A_{k,\lambda} ^{(k_1,k_2,\ldots,k_r)}(x)_{n-k,\lambda}.
\end{equation*}
\end{theorem}

Assume that $m \in {\mathbb{N}}$ with $m\equiv1$ $({\text{mod }}2)$. Then we have
\begin{equation}\label{26}
\begin{split}
&\sum_{n=0} ^{\infty} A_{n,\lambda} ^{(k_1,k_2,\ldots,k_r)}(x)\frac{t^n}{n!}=\frac{2r!}{\e(t)+1}Li_{k_1,k_2,\ldots,k_r}\left(1-e^{-t}\right)\e^x(t)\\ 
=&\frac{2r!}{\e^m(t)+1}Li_{k_1,k_2,\ldots,k_r}\left(1-e^{-t}\right)\sum_{l=0} ^{m-1}(-1)^l\e^{l+x}(t)\\
=&\frac{r!}{t^r}Li_{k_1,k_2,\ldots,k_r}\left(1-e^{-t}\right)\sum_{l=0} ^{m-1}(-1)^l\frac{2t^r}{e_{\lambda/m}(mt)+1}e_{\lambda/m} ^{\frac{l+x}{m}}(mt)\\
=&\frac{r!}{t^r}\sum_{k=r} ^{\infty} S_1 ^{(k_1,k_2,\ldots,k_r)}(k,r)\sum_{j=k} ^{\infty} S_2(j,k)(-1)^{j-k}\frac{t^j}{j!}\sum_{l=0} ^{m-1}(-1)^l\sum_{i=0} ^{\infty}A_{i,\lambda/m} ^{(r)}\left(\frac{l+x}{m}\right)m^{i-r}\frac{t^i}{i!} \\
=&\frac{r!}{t^r}\sum_{j=r} ^{\infty} \left(\sum_{k=r} ^j S_1 ^{(k_1,k_2,\ldots,k_r)}(k,r)(-1)^{j-k}S_2(j,k)\right)\frac{t^j}{j!}\sum_{i=0} ^{\infty}\left(\sum_{l=0} ^{m-1} (-1)^lA_{i,\lambda/m} ^{(r)}\left(\frac{l+x}{m}\right)m^{i-r}\right)\frac{t^i}{i!}\\
=&\sum_{j=0} ^{\infty}\left(\sum_{k=r} ^{j+r} S_1 ^{(k_1,k_2,\ldots,k_r)}(k,r)(-1)^{j+r-k}S_2(j+r,k)\frac{r!j!}{(j+r)!}\right)\frac{t^j}{j!}\\
&\times \sum_{i=0} ^{\infty} \left(\sum_{l=0} ^{m-1}(-1)^lA_{i,\lambda/m} ^{(r)} \left(\frac{l+x}{m}\right)m^{i-r}\right)\frac{t^i}{i!}\\
=&\sum_{n=0} ^{\infty}\left(\sum_{j=0} ^n \binom{n}{j}\sum_{k=r} ^{j+r}\sum_{l=0} ^{m-1} S_1 ^{(k_1,k_2,\ldots,k_r)}(k,r)(-1)^{j+r-k}\right.\\
&\times \left. \frac{S_2(j+r,k)}{\binom{j+r}{r}}(-1)^lA_{n-j,\lambda/m} ^{(r)}\left(\frac{l+x}{m}\right)m^{n-j-r}\right)\frac{t^n}{n!}.
\end{split}
\end{equation}

Therefore, by comparing the coefficients on both sides of \eqref{26}, we obtain the following theorem.
\begin{theorem}
For $m \in {\mathbb{N}}$ with $m\equiv 1$ ({\rm{mod}} $2$), we have
\begin{equation*}
\begin{split}
A_{n,\lambda} ^{(k_1,k_2,\ldots,k_r)}(x)=&\sum_{j=0} ^n \binom{n}{j}\sum_{k=r} ^{j+r}\sum_{l=0} ^{m-1} S_1 ^{(k_1,k_2,\ldots,k_r)}(k,r)(-1)^{j+r-k}\\
&\times \frac{S_2(j+r,k)}{\binom{j+r}{r}}(-1)^lA_{n-j,\lambda/m} ^{(r)}\left(\frac{l+x}{m}\right)m^{n-j-r}.
\end{split}
\end{equation*}
\end{theorem}

\section{Conclusion}
In addition to degenerate versions of many special numbers and polynomials, the degenerate gamma function, degenerate umbral calculus and degenerate $q$-umbral calculus are introduced and a lot of interesting results about them are found in recent years.\par
In this note, we introduced the multi-Euler-Genocchi and degenerate multi-Euler-Genocchi polynomials which are multi-versions of the generalized Euler-Genocchi and generalized degenerate Euler-Genocchi polynomials. Among other things, we expressed the former by the generalized Euler-Genocchi polynomials, the multi-Stirling numbers of the first kind and Stirling numbers of the second kind, and the latter by the generalized degenerate Euler-Genocchi polynomials, the multi-Stirling numbers of the first kind and Stirling numbers of the second kind. \par
It is one of our future projects to continue to study various degenerate versions of some special numbers and polynomials and those of certain transcendental functions, and to find their applications to physics, science and engineering as well as to mathematics.

\end{document}